\documentclass[graybox]{svmult}

\usepackage{mathptmx}       
\usepackage{helvet}         
\usepackage{courier}        
\usepackage{type1cm}        
\usepackage{makeidx}         
\usepackage{graphicx}        
\usepackage{multicol}        
\usepackage[bottom]{footmisc}

\usepackage{url}
\usepackage{amssymb}
\usepackage{amsfonts}
\usepackage{amsmath}
\usepackage{enumerate}
\usepackage{dsfont}
\usepackage{longtable}
\usepackage{comment}

\newcommand{\bthe}{\begin{theorem}}
\newcommand{\ethe}{\end{theorem}}

\newcommand{\ben}{\begin{enumerate}}
\newcommand{\een}{\end{enumerate}}

\newcommand{\beq}{\begin{equation}}
\newcommand{\eeq}{\end{equation}}

\newcommand{\ble}{\begin{lemma}}
\newcommand{\ele}{\end{lemma}}

\newcommand{\bde}{\begin{definition}}
\newcommand{\ede}{\end{definition}}

\newcommand{\bco}{\begin{corollary}}
\newcommand{\eco}{\end{corollary}}

\newcommand{\bpr}{\begin{proposition}}
\newcommand{\epr}{\end{proposition}}

\newcommand{\brem}{\begin{remark}}
\newcommand{\erem}{\end{remark}}

\newcommand{\bproof}{\begin{proof}}
\newcommand{\eproof}{\end{proof}}

\newcommand{\bexam}{\begin{example}}
\newcommand{\eexam}{\end{example}}

\newcommand{\bfi}{\begin{fig}}
\newcommand{\efi}{\end{fig}}

\newcommand{\btab}{\begin{tab}}
\newcommand{\etab}{\end{tab}}

\newcommand{\beao}{\begin{eqnarray*}}
\newcommand{\eeao}{\end{eqnarray*}\noindent}

\newcommand{\beam}{\begin{eqnarray}}
\newcommand{\eeam}{\end{eqnarray}\noindent}

\newcommand{\barr}{\begin{array}}
\newcommand{\earr}{\end{array}}

\newcommand{\bdis}{\begin{displaymath}}
\newcommand{\edis}{\end{displaymath}\noindent}

\newcommand{\benu}{\begin{enumerate}}
\newcommand{\eenu}{\end{enumerate}}

\newcommand{\bit}{\begin{itemize}}
\newcommand{\eit}{\end{itemize}}

\newcommand{\bff}{\textbf}

\def\E{{\Bbb E}}

\newcommand{\dd}{\mathrm{d}}

\newcommand{\bb}{\mathrm{b}}
\newcommand{\cc}{\mathrm{c}}

\newcommand{\Comp}{\mathrm{c}}
\newcommand{\loc}{\mathrm{loc}}

\newcommand{\Leb}{\mathrm{Leb}}

\newcommand{\bbn}{\mathbb{N}}

\newcommand{\bbr}{\mathbb{R}}

\newcommand{\bbe}{\mathbb{E}}

\newcommand{\bbp}{\mathbb{P}}
\newcommand{\bbb}{\mathbb{B}}
\newcommand{\bbf}{\mathbb{F}}

\newcommand{\bone}{\mathds{1}}

\newcommand{\cals}{{\cal S}}

\newcommand{\calf}{{\cal F}}

\newcommand{\calp}{{\cal P}}
\newcommand{\calb}{{\cal B}}

\newcommand{\calg}{{\cal G}}

\newcommand{\calv}{{\cal V}}
\newcommand{\calo}{{\cal O}}

\newcommand{\la}{{\lambda}}
\newcommand{\La}{{\Lambda}}
\newcommand{\ga}{{\gamma}}
\newcommand{\Ga}{{\Gamma}}
\newcommand{\si}{{\sigma}}

\newcommand{\vareps}{\varepsilon}

\newcommand{\Om}{{\Omega}}

\allowdisplaybreaks

\spdefaulttheorem{alg}{Algorithm}{\bfseries}{}

\makeindex             

\begin{document}

\title*{Simulation of Stochastic Volterra Equations Driven by Space--Time L\'evy Noise}
\author{Bohan Chen \and Carsten Chong \and Claudia Kl\"uppelberg}
\institute{Bohan Chen \and Carsten Chong \and Claudia Kl\"uppelberg \at Technische Universit\"at M\"unchen \\ Boltzmannstra\ss e 3, 85748 Garching, Germany
\and Bohan Chen \at \email{bohan.chen@mytum.de}
\and Carsten Chong \at \email{carsten.chong@tum.de}
\and Claudia Kl\"uppelberg \at \email{cklu@ma.tum.de}}

\maketitle

\abstract{In this paper we investigate two numerical schemes for the simulation of stochastic Volterra equations driven by space--time L\'evy noise of pure-jump type. The first one is based on truncating the small jumps of the noise, while the second one relies on series representation techniques for infinitely divisible random variables. Under reasonable assumptions, we prove for both methods $L^p$- and almost sure convergence of the approximations to the true solution of the Volterra equation. We give explicit convergence rates in terms of the Volterra kernel and the characteristics of the noise. A simulation study visualizes the most important path properties of the investigated processes. \keywords{Simulation of SPDEs, Simulation of stochastic Volterra equations, Space--time L\'evy noise, Stochastic heat equation, Stochastic partial differential equation}}

\vspace{\baselineskip}
\noindent\bff{Mathematics Subject Classifications (2010): 60H35, 65C30, 60H20, 60H15, 60G57, 60G51}

\section{Introduction}\label{s1}

The aim of this paper is to investigate different simulation techniques for stochastic Volterra equations (SVEs) of the form
\beq\label{stochint} Y(t,x) = Y_0(t,x) + \int_0^t \int_{\bbr^d} G(t,x;s,y)\si(Y(s,y))\,\Lambda(\dd s,\dd y)\;, \quad(t,x)\in \bbr_+\times\bbr^d\;,\eeq
where $G$ is a deterministic kernel function, $\si$ a Lipschitz coefficient and $\Lambda$ a L\'evy basis on $\bbr_+\times\bbr^d$ of pure-jump type with no Gaussian part. In the purely temporal case where no space is involved and the kernel $G$ is sufficiently regular on the diagonal $\{(t;s)\in \bbr_+\times\bbr_+\colon t=s\}$, the existence and uniqueness of the solution $Y$ to \eqref{stochint} are established for general semimartingale integrators in \cite{Protter85}. The space--time case \eqref{stochint} is treated in \cite{Chong14b} for quite general L\'evy bases. In particular, $G$ is allowed to be singular on the diagonal, which typically happens in the context of stochastic partial differential equations (SPDEs) where $G$ is the Green's function of the underlying differential operator. More details on the connection between SPDEs and the SVE \eqref{stochint} are presented in Sect.~\ref{s2}, or can be found in \cite{BN11-2, Chong14b, Walsh86}.


Since in most cases there exists no explicit solution formula for the SVE \eqref{stochint}, it is a natural task to develop appropriate simulation algorithms. For SPDEs driven by Gaussian noise, research on this topic is rather far advanced, see e.g. \cite{Davie00, Gyoengy99, Walsh05}. However, for SPDEs driven by jump noises such as non-Gaussian L\'evy bases, the related literature is considerably smaller, see \cite{Barth12} and the work of Hausenblas and coauthors \cite{Dunst12, Hausenblas08, Hausenblas06}. The case $\si\equiv1$ has been treated in \cite{Chen14}. The contribution of our paper can be summarized as follows:
\bit
	\item We propose and analyze two approximation schemes for \eqref{stochint}, each of which replaces the original noise by a truncated noise that only has finitely many atoms on compact subsets of $\bbr_+\times\bbr^d$. For the first scheme, we simply cut off all jumps whose size is smaller than a constant. For the second scheme, we use series representation techniques for the noise as in \cite{Rosinski90} such that the jumps to be dropped off are chosen randomly. Both methods have already been applied successfully to the simulation of L\'evy processes, cf. \cite{Asmussen01,Rosinski01}.
	\item In the case where $G$ originates from an SPDE, the crucial difference of our numerical schemes to the Euler or finite element methods in the references mentioned before is that we do not simulate small space--time increments of the noise but successively the true jumps of the L\'evy basis, which is an easier task given that one usually only knows the underlying L\'evy measure. 
It is important to recognize that this is only possible because the noise $\La$ is of pure-jump type, and contains neither a Gaussian part nor a drift. We shall point out in Sect.~\ref{s6} how to relax this assumption.
\eit

The remaining article is organized as follows: Section~2 gives the necessary background for the SVE \eqref{stochint}. In particular, we present sufficient conditions for the existence and uniqueness of solutions, and address the connection between \eqref{stochint} and SPDEs. In Sect.~3 we construct approximations to the solution $Y$ of \eqref{stochint} by truncating the small jumps of the L\'evy basis. We prove in Thm.~\ref{conv1} their $L^p$-convergence, and in some cases also their almost sure (a.s.) convergence to the target process $Y$. In Sect.~4 we approximate the driving L\'evy basis using series representation methods. This leads to an algorithm that produces approximations again converging in the $L^p$-sense, sometimes also almost surely, to $Y$, see Thm.~\ref{conv2}. In both theorems, we find explicit $L^p$-convergence rates that only depend on the kernel $G$ and the characteristics of $\La$. Section~5 presents a simulation study for the stochastic heat equation which highlights the typical path behaviour of stochastic Volterra equations. The final Sect.~6 compares the two simulation algorithms developed in this paper and discusses some further directions of the topic.

\section{Preliminaries}\label{s2}
We start with a summary of notations that will be employed in this paper. 
\begin{longtable}{p{1.5 cm} p{10 cm}}
$\bbr_+$ & the set $[0,\infty)$ of \emph{positive} real numbers\\
$\bbn$ & the natural numbers $\{1, 2, \ldots\}$\\
$\bbb$ & a stochastic basis $(\Om,\calf,\bbf=(\calf_t)_{t\in\bbr_+},\bbp)$ satisfying the usual hypotheses of completeness and right-continuity\\
$\bar\Om$,$\tilde\Om$ & $\bar\Om:=\Om\times\bbr_+$ and $\tilde\Om:=\Om\times\bbr_+\times\bbr^d$ where $d\in\bbn$ \\
$\calb(\bbr^d)$ & the Borel $\si$-field on $\bbr^d$\\
$\tilde\calb_\bb$ & the collection of all bounded Borel sets of $\bbr_+\times\bbr^d$ \\
$\calp$ & the predictable $\si$-field on $\bbb$ \emph{or} the collection of all predictable processes $\bar\Om\to\bbr$\\
$\tilde\calp$ & the product $\calp\otimes\calb(\bbr^d)$ \emph{or} the collection all $\calp\otimes\calb(\bbr^d)$-measurable processes $\tilde\Om\to\bbr$\\
$\tilde\calp_\bb$ & the collection of sets in $\tilde\calp$ which are a subset of $\Om\times[0,k]\times[-k,k]^d$ for some $k\in\bbn$ \\
$p^\ast$ & $p\vee1$\\
$L^p$ & the space $L^p(\Om,\calf,\bbp)$, $p\in(0,\infty]$, endowed with the topology induced by $\|X\|_{L^p}:=\bbe[|X|^p]^{1/{p^\ast}}$\\
$L^0$ & the space $L^0(\Om,\calf,\bbp)$ of all random variables on $\bbb$ endowed with the topology of convergence in probability\\
$B^p_\loc$ & the set of all $Y\in\tilde\calp$ for which $\|Y(t,x)\|_{L^p}$ is uniformly bounded on $[0,T]\times\bbr^d$ for all $T\in\bbr_+$ ($p\in(0,\infty]$)\\
$A^\Comp$ & the complement of $A$ within the superset it belongs to (which will be clear from the context)\\
$A-B$ & $\{x-y\colon x\in A, y\in B\}$\\
$-A$ & $\{-x\colon x\in A\}$\\
$\Leb$ & the Lebesgue measure on $\bbr^d$ ($d$ should be clear from the context)\\
$\|\cdot\|$ & the Euclidean norm on $\bbr^d$.\\
$C, C(T)$ & two generic constants in $\bbr_+$, one dependent and one independent of $T$, whose values we do not care of and may therefore change from one place to the other
\end{longtable}

We suppose that the stochastic basis $\bbb$ supports a \emph{L\'evy basis}, that is, a mapping $\La\colon \tilde\calp_\bb \to L^0$ with the following properties:
\begin{itemize}
  \item $\La(\emptyset)=0$ a.s.
  \item For all pairwise disjoint sets $(A_i)_{i\in\mathbb{N}}\subset\tilde\calp_\bb$ with $\bigcup_{i=1}^{\infty} A_i\in\tilde\calp_\bb$ we have
      \beq
      \La\left(\bigcup_{i=1}^{\infty}A_i\right)=\sum_{i=1}^{\infty}\La(A_i) \quad \text{in } L^0\;.\label{eq:1}
      \eeq
  \item $(\La(\Om\times B_i))_{i\in\mathbb{N}}$ is a sequence of independent random variables if $(B_i)_{i\in\bbn}$ are pairwise disjoint sets in $\tilde\calb_\bb$.
  \item For every $B\in\mathcal{B}_\bb$, $\La(\Om\times B)$ has an infinitely divisible distribution.
  \item $\La(A)$ is $\calf_t$-measurable when $A\in\tilde\calp_\bb$ and $A\subseteq \Om\times[0,t]\times\bbr^d$ for $t\in\bbr_+$.
  \item For every $t\in\bbr_+$, $A\in\tilde\calp_\bb$ and $\Om_0\in\calf_t$ we have a.s.
  \[ \La(A\cap(\Om_0\times(t,\infty)\times\bbr^d))=\bone_{\Om_0} \La(A\cap(\Om\times(t,\infty)\times\bbr^d))\;. \]
\end{itemize}

Just as L\'evy processes are semimartingales and thus allow for an It\^o integration theory, L\'evy bases belong to the class of $L^0$-valued $\si$-finite random measures. Therefore, it is possible to define the stochastic integral
	\[ \int_{\bbr_+\times\bbr^d} H(s,y)\,\La(\dd s,\dd y) \]
for $H\in\tilde\calp$ that are \emph{integrable} with respect to $\La$, see \cite{Chong14} for the details.

Similarly to L\'evy processes, there exist two notions of characteristics for L\'evy bases: one going back to \cite[Prop.~2.1]{Rajput89} that is based on the L\'evy-Khintchine formula and is independent of $\bbf$, and a filtration-based one that is useful for stochastic analysis \cite[Thm.~3.2]{Chong14}. For the whole paper, we will assume that both notions coincide such that $\La$ has a canonical decomposition under the filtration $\bbf$ of the form
	\begin{align*} \La(\dd t,\dd x) &= B(\dd t,\dd x) + \La^\cc(\dd t,\dd x) + \int_\bbr z\bone_{\{|z|\leq 1\}} \,(\mu-\nu)(\dd t,\dd x,\dd z)\\
	&\quad + \int_\bbr z\bone_{\{|z|>1\}}\,\mu(\dd t,\dd x,\dd z)\;, \end{align*}
where $B$ is a $\si$-finite signed Borel measure on $\bbr_+\times\bbr^d$, $\La^\cc$ a L\'evy basis such that $\La(\Om\times B)$ is normally distributed with mean $0$ and variance $C(B)$ for all $B\in\tilde\calb_\bb$, and $\mu$ a Poisson measure on $\bbr_+\times\bbr^d$ relative to $\bbf$ with intensity measure $\nu$ (cf. \cite[Def.~II.1.20]{Jacod03}). There exists also a $\si$-finite Borel measure $\la$ on $\bbr_+\times\bbr^d$ such that
	\begin{align} B(\dd t,\dd x) &= b(t,x)\,\la(\dd t,\dd x)\;,\quad C(\dd t,\dd x)=c(t,x)\,\la(\dd t,\dd x)\quad\text{and}\nonumber\\
	\nu(\dd t,\dd x,\dd z)&=\pi(t,x,\dd z)\,\la(\dd t,\dd x) \label{charLa}\end{align}
	with two functions $b\colon \bbr_+\times\bbr^d \to \bbr$ and $c\colon \bbr_+\times\bbr^d\to \bbr_+$ as well as a transition kernel $\pi$ from $(\bbr_+\times\bbr^d,\calb(\bbr_+\times\bbr^d))$ to $(\bbr,\calb(\bbr))$ such that $\pi(t,x,\cdot)$ is a L\'evy measure for each $(t,x)\in\bbr_+\times\bbr^d$. 
		

We have already mentioned in the introduction that we will assume
\beq\label{Gauss0} C = 0 \eeq
throughout the paper. For simplicity we will also make two further assumptions: first, that there exist $b\in\bbr$ and a L\'evy measure $\pi$ such that for all $(t,x)\in\bbr_+\times\bbr^d$ we have
\beq\label{invar} b(t,x)=b\;, \quad \pi(t,x,\cdot)=\pi\quad\text{and}\quad \la(\dd t,\dd x)=\dd(t,x)\; ;\eeq
second, that 
\beq\label{finvarsymm} \La \in \cals \cup \calv_0\;, \eeq
where $\cals$ is the collection of all \emph{symmetric} L\'evy bases and $\calv_0$ is the class of L\'evy bases with \emph{locally finite variation and no drift}, defined by the property that
\[ \int_\bbr |z|\bone_{\{|z|\leq1\}} \,\pi(\dd z)<\infty\;,\quad\text{and}\quad b_0:=b-\int_\bbr z\bone_{\{|z|\leq1\}} \,\pi(\dd z)=0\;.\]
Furthermore, if $\pi$ has a finite first moment, that is, 
\beq\label{meanfin} \int_\bbr |z|\bone_{\{|z|>1\}} \,\pi(\dd z)<\infty\;, \eeq
we define 
	\begin{align*} B_1(\dd t,\dd x)&:=b_1\,\dd(t,x)\;,\quad b_1:=b+\int_\bbr z\bone_{\{|z|>1\}} \,\pi(\dd z)\;,\\
	M(\dd t,\dd x) &:= \La(\dd t,\dd x) - B_1(\dd t,\dd x) = \int_\bbr z\,(\mu-\nu)(\dd t,\dd x,\dd z)\;. \end{align*} 

Next, let us summarize the most important facts regarding the SVE \eqref{stochint}. All details that are not explained can be found in \cite{Chong14b}. 
First, many SPDEs of evolution type driven by L\'evy noise can be written in terms of \eqref{stochint}, where $G$ is the Green's function of the corresponding differential operator. Most prominently, taking $G$ being the heat kernel in $\bbr^d$, \eqref{stochint} is the so-called \emph{mild formulation} of the stochastic heat equation (with constant coefficients and multiplicative noise). Typically for parabolic equations, the heat kernel is very smooth in general but explodes on the diagonal $t=s$ and $x=y$. In fact, it is only $p$-fold integrable on $[0,T]\times\bbr^d$ for $p<1+2/d$. In particular, as soon as $d\geq2$, it is not square-integrable, and as a consequence, no solution to the stochastic heat equation in the form \eqref{stochint} will exist for L\'evy noises with non-zero Gaussian component. This is another reason for including assumption \eqref{Gauss0} in this paper.

Second, let us address the existence and uniqueness problem for \eqref{stochint}. By a \emph{solution} to this equation we mean a predictable process $Y\in\tilde\calp$ such that for all $(t,x)\in\bbr_+\times\bbr^d$, the stochastic integral on the right-hand side of \eqref{stochint} is well defined and the equation itself for each $(t,x)\in[0,T]\times\bbr^d$ holds a.s. We identify two solutions as soon as they are modifications of each other. Given a number $p\in(0,2]$, the following conditions guarantee a unique solution to \eqref{stochint} in $B^p_\loc$ by \cite[Thm.~3.1]{Chong14b}:
\benu
	\item[A1.] $Y_0\in B^p_\loc$ is independent of $\La$.
	\item[A2.] $\si\colon \bbr\to\bbr$ is Lipschitz continuous, that is, there exists $C\in\bbr_+$ such that 
	\beq\label{Lip}|\si(x)-\si(y)|\leq C|x-y|\;,\quad x,y\in\bbr\;. \eeq
	\item[A3.] $G\colon (\bbr_+\times\bbr^d)^2\to\bbr$ is a measurable function with $G(t,\cdot;s,\cdot)\equiv0$ for $s>t$.
	\item[A4.] $\La$ satisfies \eqref{charLa}--\eqref{finvarsymm} and
	\beq\label{intp} \int_{\bbr} |z|^p\,\pi(\dd z)<\infty\;. \eeq
	\item[A5.] If we define for $(t,x),(s,y)\in\bbr_+\times\bbr^d$
	\beq\label{tildeG} \tilde G(t,x;s,y):= |G(t,x;s,y)|\bone_{\{p>1,\La\notin\cals\}} + |G(t,x;s,y)|^p\;,\eeq 
	then we have for all $T\in\bbr_+$
	\beq\label{integrab} \sup_{(t,x)\in[0,T]\times\bbr^d} \int_0^T \int_{\bbr^d} \tilde G(t,x;s,y)\,\dd(s,y) <\infty\;. \eeq
	\item[A6.] For all $\vareps>0$ and $T\in\bbr_+$ there exist $k\in\bbn$ and a partition $0=t_0<\ldots<t_k=T$ such that
	\beq\label{subdiv} \sup_{(t,x)\in[0,T]\times\bbr^d} \sup_{i=1,\ldots,k} \int_{t_{i-1}}^{t_i} \int_{\bbr^d} \tilde G(t,x;s,y) \,\dd(s, y) < \vareps\;. \eeq
\eenu
Apart from A1--A6, we will add another assumptions in this paper:
\benu
	\item[A7.] There exists a sequence $(U^N)_{N\in\bbn}$ of compact sets increasing to $\bbr^d$ such that for all $T\in\bbr_+$ and compact sets $K\subseteq\bbr^d$ we have, as $N\to\infty$,
	\begin{align}  r^N_1(T,K)&:=\sup_{(t,x)\in[0,T]\times K} \Bigg(\int_0^t\int_{(U^N)^\Comp} |G(t,x;s,y)|\bone_{\{p>1,\La\notin\cals\}}\,\dd (s, y)\nonumber\\
 &\quad+ \left(\int_0^t\int_{(U^N)^\Comp} |G(t,x;s,y)|^p\,\dd(s,y)\right)^{1/{p^\ast}}\Bigg)\to 0\;.\label{rN1}
	\end{align} 
\eenu

Conditions A6 and A7 are automatically satisfied if $|G(t,x;s,y)|\leq g(t-s,x-y)$ for some measurable function $g$ and A5 holds with $G$ replaced by $g$. For A6 see \cite[Rem.~3.3(3)]{Chong14b}; for A7 choose $U^N:=\{x\in\bbr^d\colon \|x\|\leq N\}$ such that for $p>1$
\begin{align*} &~\sup_{(t,x)\in[0,T]\times K} \int_0^t\int_{(U^N)^\Comp} |G(t,x;s,y)|^p\,\dd (s,y)\\
\leq&~ \sup_{(t,x)\in[0,T]\times K} \int_0^t\int_{(U^N)^\Comp} g^p(t-s,x-y)\,\dd (s,y)\\
 =&~\sup_{x\in K} \int_0^T \int_{x-(U^N)^\Comp} g^p(s,y)\,\dd(s,y)\leq \int_0^T\int_{K-(U^N)^\Comp} g^p(s,y)\,\dd(s,y)\\
 \to&~0\quad\text{as}\quad N\to\infty  \end{align*}
 by the fact that $K-(U^N)^\Comp\downarrow0$. A similar calculation applies to the case $p\in(0,1]$ and the first term in $r^N_1(T,K)$. 
\bexam We conclude this section with the stochastic heat equation in $\bbr^d$, whose mild formulation is given by the SVE \eqref{stochint} with
\beq\label{hk} G(t,x;s,y)=g(t-s,x-y)\;,\quad g(t,x)=\frac{\exp(-\|x\|^2/(4t))}{(4\pi t)^{d/2}}\bone_{[0,t)}(s) \eeq
for $(t,x),(s,y)\in\bbr_+\times\bbr^d$. We assume that $Y_0$ and $\si$ satisfy conditions A1 and A2, respectively. Furthermore, we suppose that \eqref{charLa}--\eqref{finvarsymm} are valid, and that \eqref{intp} holds with some $p\in(0,1+2/d)$. It is straightforward to show that then A3--A6 are satisfied with the same $p$. Let us estimate the rate $r^N_1(T,K)$ for $T\in\bbr_+$, $K:=\{\|x\|\leq R\}$ with $R\in\bbn$, and $U^N:=\{\|x\|\leq N\}$. We first consider the case $p\leq1$ or $\La\in\cals$. Since $K-(U^N)^\Comp=(U^{N-R})^\Comp$ for $N\geq R$, the calculations after A7 yield ($\Gamma(\cdot,\cdot)$ denotes the upper incomplete gamma function and $p(d):=1+(1-p)d/2$)
\begin{align}
	(r^N_1(T,K))^{p^\ast}&\leq \int_0^T \int_{(U^{N-R})^\Comp} g^p(t,x)\,\dd(t,x)=\int_0^T\int_{N-R}^\infty \frac{\exp(-pr^2/(4t))}{(4\pi t)^{pd/2}} r^{d-1}\,\dd r\,\dd t\nonumber\\
	&=C\int_0^T t^{p(d)-1}\Ga\left(\frac{d}{2},\frac{p(N-R)^2}{4t}\right)\,\dd t\nonumber\\
	&=C(T) \Bigg((p(d))^{-1}\Ga\left(\frac{d}{2},\frac{p(N-R)^2}{4T}\right)-\left(\frac{p(N-R)^2}{4T}\right)^{p(d)}\nonumber\\
	&\quad\times\Ga\left(\frac{d}{2}-p(d),\frac{p(N-R)^2}{4T}\right)\Bigg)\nonumber\\
	&\leq C(T) \exp\left(-\frac{p(N-R)^2}{4T}\right)(N-R)^{d-2}\;, \label{rN1-ex} 
\end{align}
which tends to $0$ exponentially fast as $N\to\infty$. If $p>1$ and $\La\notin\cals$, it follows from formula \eqref{rN1} that we need an extra summand for $r^N_1(T,K)$, namely \eqref{rN1-ex} with $p=1$.
\eexam

\section{Truncation of Small Jumps}\label{s3}

In this section we approximate equation \eqref{stochint} by cutting off the small jumps of $\La$. To this end, we first define for each $N\in\bbn$ \beq\label{GN} G^N(t,x;s,y):=G(t,x;s,y)\bone_{U^N}(y)\;,\quad (t,x),(s,y)\in\bbr_+\times\bbr^d\;,\eeq
where the meaning of the sets $U^N$ is explained in A7. Furthermore, we introduce
	\beq	\label{rNs}  r^N_2:= \left(\int_{[-\vareps^N,\vareps^N]} |z|^p\,\pi(\dd z)\right)^{1/{p^\ast}}\;,\quad r^N_3:=\left|\int_{[-\vareps^N,\vareps^N]} z\bone_{\{p>1,\La\notin\cals\}}\,\pi(\dd z)\right|\;, \eeq
		where $(\vareps^N)_{N\in\bbn}\subseteq(0,1)$ satisfies $\vareps^N\to0$ as $N\to\infty$. Condition A4 implies that $r^N_2,r^N_3\to0$ as $N\to\infty$. Next, defining truncations of the L\'evy basis $\La$ by
\beq \La^N(\dd t,\dd x) 
	 := \int_{[-\vareps^N,\vareps^N]^\Comp} z\,\mu(\dd t,\dd x,\dd z)\;, \label{LaN} \eeq
our approximation scheme for the solution $Y$ to \eqref{stochint} is given as:
\beq\label{scheme1} Y^N(t,x) := Y_0(t,x) + \int_0^t \int_{\bbr^d} G^N(t,x;s,y)\si(Y^N(s,y))\,\La^N(\dd s,\dd y) \eeq
for $(t,x)\in \bbr_+\times\bbr^d$.
Indeed, $Y^N$ can be simulated exactly because for all $T\in\bbr_+$ the truncation $\La^N$ only has a finite (random) number $R^N(T)$ of jumps on $[0,T]\times U^N$, say at the space--time locations $(\tau^N_i,\xi^N_i)$ with sizes $J^N_i$. This implies that we have the following alternative representation of $Y^N(t,x)$ for $(t,x)\in[0,T]\times\bbr^d$:
\beq\label{scheme1var} Y^N(t,x) = Y_0(t,x) + \sum_{i=1}^{R^N(T)} G(t,x;\tau_i^N,\xi_i^N)\si(Y^N(\tau_i^N,\xi_i^N))J^N_i\bone_{\{\tau^N_i<t\}}\;. \eeq
What remains to do is to simulate $Y^N(\tau_i^N,\xi_i^N)$, $i=1,\ldots, R^N(T)$, iteratively, from which the values $Y(t,x)$ for all other $(t,x)\in[0,T]\times\bbr^d$ can be computed. 


The following algorithm summarizes up the simulation procedure:
\begin{alg}\label{alg1} Consider a finite grid $\calg$ that is a subset of $[0,T]\times\bbr^d$. For each step $N$ proceed as follows:
\begin{enumerate}
  \item Draw a Poisson random variable $R^N(T)$ with intensity 
  \[ R^N(T) := \int_0^T\int_{U^N}\int_\bbr \bone_{\{z\in[-\vareps^N,\vareps^N]^\Comp\}}\,\nu(\dd t,\dd x,\dd z) = T\Leb(U^N)\pi\big([-\vareps^N,\vareps^N]^\Comp\big)\;.\]
  \item For $i=1,\ldots,R^N(T)$:
  	\benu
  		\item Draw a pair $(\tau^N_i,\xi^N_i)$ with uniform distribution from $[0,T]\times U^N$.
  		\item Draw $J^N_i$ from $[-\vareps^N,\vareps^N]^\Comp$ with distribution $\pi/\pi\big([-\vareps^N,\vareps^N]^\Comp\big)$.
  	\eenu
  \item For each $i=1,\ldots,R^N(T)$ and $(t,x)\in\calg$ simulate $Y_0(\tau^N_i,\xi^N_i)$ and $Y_0(t,x)$.
  \item For each $i=1,\ldots,R^N(T)$ set 
  \[ Y^N(\tau^N_i,\xi^N_i) := Y_0(\tau^N_i,\xi^N_i) + \sum_{j=1}^{i-1} G(\tau^N_i,\xi^N_i;\tau^N_j,\xi^N_j)\si(Y^N(\tau^N_j,\xi^N_j)) J^N_j\;. \]
  \item For each $(t,x)\in\calg$ define $Y^N(t,x)$ via \eqref{scheme1var}.  
\end{enumerate}
\end{alg}

The next theorem determines the convergence behaviour of the scheme \eqref{scheme1} to the true solution $Y$ to \eqref{stochint}. 
\bthe\label{conv1} Grant assumptions A1--A7 under which the SVE \eqref{stochint} has a unique solution in $B^p_\loc$. Then $Y^N$ as defined in \eqref{scheme1} belongs to $B^p_\loc$ for all $N\in\bbn$, and for all $T\in\bbr_+$ and compact sets $K\subseteq\bbr^d$ there exists a constant $C(T)\in\bbr_+$ independent of $N$ and $K$ such that
\beq\label{convtrunc} \sup_{(t,x)\in[0,T]\times K} \|Y(t,x)-Y^N(t,x)\|_{L^p} \leq C(T) (r^N_1(T,K) + r^N_2 + r^N_3)\;. \eeq
Furthermore, if $\sum_{n=1}^\infty (r^N_1(T,K)+r^N_2 +r^N_3)^{p^\ast} < \infty$ is fulfilled, then we also have for all $(t,x)\in[0,T]\times K$ that $Y^N(t,x)\to Y(t,x)$ a.s. as $N\to\infty$.
\ethe


\bproof
It is obvious that $|G^N|\leq|G|$ pointwise and that we have $\nu^N\leq \nu$ for the third characteristic $\nu^N$ of $\La^N$. Thus, A1--A6 are still satisfied when $G$ and $\nu$ are replaced by $G^N$ and $\nu^N$ (if $\La\in\cals$, also $\La^N\in\cals$). So $Y^N$ as a solution to \eqref{stochint} with $G^N$ and $\La^N$ instead of $G$ and $\La$ belongs to $B^p_\loc$ as well. Moreover, for all $T\in\bbr_+$ there exists $C(T)\in\bbr_+$ independent of $N\in\bbn$ such that 
\beq\label{YNY} \sup_{(t,x)\in[0,T]\times\bbr^d} \|Y^N(t,x)\|_{L^p} \leq C(T)\;,\quad N\in\bbn\;. \eeq
We only sketch the proof for this statement. In fact, using \cite[Lem.~6.1(1)]{Chong14b} it can be shown that the left-hand side of \eqref{YNY} satisfies an inequality of the same type as in Lem.~6.4(3) of the same paper. In particular, it is bounded by a constant $C^N(T)$ that depends on $N$ only through $|G^N|$ and $\nu^N$, and that this constant is only increased if we replace $|G^N|$ and $\nu^N$ by the larger $|G|$ and $\nu$. In this way, we obtain an upper bound $C(T)$ that does not depend on $N$.

Next, we prove the convergence of $Y^N$ to $Y$ as stated in \eqref{convtrunc}. We have
\begin{align}
 Y(t,x)-Y^N(t,x) &= \int_0^t\int_{\bbr^d} [G(t,x;s,y)-G^N(t,x;s,y)] \si(Y(s,y))\,\La(\dd s,\dd y) \nonumber \\
 &\quad+ \int_0^t\int_{\bbr^d} G^N(t,x;s,y)[\si(Y(s,y))-\si(Y^N(s,y))]\,\La(\dd s,\dd y) \nonumber\\
 &\quad + \int_0^t\int_{\bbr^d} G^N(t,x;s,y)\si(Y^N(s,y))\,(\La-\La^N)(\dd s,\dd y) \nonumber\\
 &=: I^N_1(t,x) + I^N_2(t,x) + I^N_3(t,x)\;, \quad (t,x)\in\bbr_+\times\bbr^d\;. \label{INs}
\end{align}
If $p>1$, we have by \eqref{Lip}, H\"older's inequality and the Burkholder-Davis-Gundy-inequality
\begin{align}
 \|I^N_2(t,x)\|_{L^p} &\leq  \left\|\int_0^t\int_{\bbr^d} G^N(t,x;s,y) [\si(Y(s,y))-\si(Y^N(s,y))]\,B_1(\dd s,\dd y)\right\|_{L^p} \nonumber\\
 &\quad+ \left\|\int_0^t\int_{\bbr^d} G^N(t,x;s,y) [\si(Y(s,y))-\si(Y^N(s,y))]\,M(\dd s,\dd y)\right\|_{L^p}\nonumber\\
 &\leq C \Bigg(\left(\int_0^t\int_{\bbr^d} |G(t,x;s,y)|\,|B_1|(\dd s,\dd y)\right)^{p-1}\nonumber\\
 &\quad\times\int_0^t\int_{\bbr^d} |G(t,x;s,y)|\|Y(s,y)-Y^N(s,y)\|^p_{L^p}\,|B_1|(\dd s,\dd y)\Bigg)^{1/p}\nonumber\\
 &\quad+C \left(\int_0^t\int_{\bbr^d} |G(t,x;s,y)|^p\|Y(s,y)-Y^N(s,y)\|_{L^p}^p\,\dd(s,y) \right)^{1/p}\;. \label{IN2-1}
\end{align}
If $p\in(0,1]$, we have $\La\in\calv_0$ by \eqref{finvarsymm} and \eqref{intp}, and thus Jensen's inequality gives
\begin{align}
 \|I^N_2(t,x)\|_{L^p} &=\bbe\left[\left(\int_0^t\int_{\bbr^d} G^N(t,x;s,y) [\si(Y(s,y))-\si(Y^N(s,y))]z\,\mu(\dd s,\dd y,\dd z)\right)^p\right]\nonumber\\
 &\leq \bbe\left[\int_0^t\int_{\bbr^d} |G^N(t,x;s,y) [\si(Y(s,y))-\si(Y^N(s,y))]z|^p\,\nu(\dd s,\dd y,\dd z)\right]\nonumber\\
 &\leq C \int_0^t\int_{\bbr^d} |G(t,x;s,y)|^p\|Y(s,y)-Y^N(s,y)\|_{L^p}\,\dd(s,y)\;.\label{IN2-2}
\end{align}
Inserting \eqref{IN2-1} and \eqref{IN2-2} back into \eqref{INs}, we have for $v^N(t,x):=\|Y(t,x)-Y^N(t,x)\|_{L^p}$
\begin{align*}
v^N(t,x) &\leq C(T)\Bigg(\left(\int_0^t\int_{\bbr^d} |G(t,x;s,y)|\bone_{\{p>1,\La\notin\cals\}} (v^N(s,y))^p\,\dd(s,y)\right)^{1/p} \\
&\quad + \left(\int_0^t\int_{\bbr^d}|G(t,x;s,y)|^p (v^N(s,y))^{p^\ast}\,\dd(s,y)\right)^{1/{p^\ast}} \Bigg)\\
&\quad + \|I^N_1(t,x) + I^N_3(t,x)\|_{L^p}\;,\quad (t,x)\in[0,T]\times\bbr^d\;.
\end{align*}
By a Gronwall-type estimate, which is possible because of A5 (see the proof of \cite[Thm.~4.7(3)]{Chong14b} for an elaboration of an argument of this type), we conclude
\[ \sup_{(t,x)\in[0,T]\times K} v^N(t,x) \leq C(T) \sup_{(t,x)\in[0,T]\times K} \|I^N_1(t,x)+I^N_3(t,x)\|_{L^p}\;.  \]
where $C(T)$ does not depend on $K$ because of \eqref{integrab}. For $I^N_1(t,x)$ we have for $p>1$
\begin{align}
 \|I^N_1(t,x)\|_{L^p} &\leq  \left\|\int_0^t\int_{\bbr^d} [G(t,x;s,y)-G^N(t,x;s,y)] \si(Y(s,y))\,B_1(\dd s,\dd y)\right\|_{L^p} \nonumber\\
 &\quad+ \left\|\int_0^t\int_{\bbr^d} [G(t,x;s,y)-G^N(t,x;s,y)] \si(Y(s,y))\,M(\dd s,\dd y)\right\|_{L^p}\nonumber\\
 &\leq C\Bigg(1+\sup_{(t,x)\in[0,T]\times\bbr^d} \|Y(t,x)\|_{L^p}\Bigg) \Bigg(\int_0^t\int_{(U^N)^\Comp} |G(t,x;s,y)|\,|B_1|(\dd s,\dd y)\nonumber\\
 &\quad+ \left(\int_0^t\int_{(U^N)^\Comp} |G(t,x;s,y)z|^p\,\nu(\dd s,\dd y,\dd z) \right)^{1/p}\Bigg)\nonumber\\
 &\leq C(T) r^N_1(T,K)\;, \label{IN1}
\end{align}
uniformly in $(t,x)\in[0,T]\times K$. In similar fashion one proves the estimate \eqref{IN1} for $p\in(0,1]$, perhaps with a different $C(T)$. Next, when $p>1$, \eqref{YNY} implies
\begin{align*}
 \|I^N_3(t,x)\|_{L^p} &=  \left\|\int_0^t\int_{\bbr^d}\int_{[-\vareps^N,\vareps^N]} G^N(t,x;s,y) \si(Y^N(s,y))z\,(\mu-\nu)(\dd s,\dd y,\dd z)\right\|_{L^p}\nonumber \\
 &\quad+\left\|\int_0^t\int_{\bbr^d} \int_{[-\vareps^N,\vareps^N]} G^N(t,x;s,y) \si(Y^N(s,y)) z\bone_{\{\La\notin\cals\}} \,\nu(\dd s,\dd y,\dd z)\right\|_{L^p}\\
 &\leq C(T)\Bigg(\left(\int_0^t\int_{\bbr^d} \int_{[-\vareps^N,\vareps^N]} |G(t,x;s,y)z|^p \,\pi(\dd z)\,\dd(s,y)\right)^{1/p} \\
 &\quad+ \left|\int_{[-\vareps^N,\vareps^N]} z\bone_{\{\La\notin\cals\}} \,\pi(\dd z)\right|\int_0^t\int_{\bbr^d} |G(t,x;s,y)|\bone_{\{\La\notin\cals\}}\,\dd(s,y)\Bigg)\\
 &\leq C(T) (r^N_2+r^N_3)\;.
\end{align*}
The case $p\in(0,1]$ can be treated similarly, cf. the estimation of $I^N_2(t,x)$ above.

It remains to prove that for each $(t,x)\in[0,T]\times K$ the convergence of $Y^N(t,x)$ to $Y(t,x)$ is almost sure when $r^N_1(T,K)$, $r^N_2$ and $r^N_3$ are ${p^\ast}$-summable. To this end, choose an arbitrary sequence $(a_N)_{N\in\bbn}\subseteq(0,1)$ converging to $0$ such that
\[  \sum_{N=1}^\infty A_N <\infty\quad\text{with}\quad A_N:=\frac{(r^N_1(T,K)+r^N_2+r^N_3)^{p^\ast}}{a^p_N}\;. \]
Such a sequence always exists, see \cite[Thm.~175.4]{Knopp90}, for example. So by \eqref{convtrunc} and Chebyshev's inequality we derive
\[ \bbp\left[ |Y(t,x)-Y^N(t,x)| \geq a_N\right] \leq \frac{\|Y(t,x)-Y^N(t,x)\|^{p^\ast}_{L^p}}{a_N^p} \leq C(T)A_N\;.  \]
Our assertion now follows from the Borel-Cantelli lemma. \qed
\eproof

\bexam\label{ex1} The rates $r^N_2$ and $r^N_3$ from \eqref{rNs} only depend on the underlying L\'evy measure $\pi$. Let $p,q\in(0,2]$ with $q<p$, and assume that $\int_{[-1,1]} |z|^q \,\pi(\dd z)<\infty$. If $\La\in\calv_0$, assume that $q<1$. Then
\begin{align*} r^N_2 &= \left(\int_{[-\vareps^N,\vareps^N]} |z|^p\,\pi(\dd z)\right)^{1/p^\ast} \leq \left(\int_{[-1,1]} |z|^q \,\nu(\dd z) (\vareps^N)^{p-q}\right)^{1/p^\ast}\\
&= \calo\left((\vareps^N)^{(p-q)/p^\ast}\right)\;,\\
r^N_3 &= \left|\int_{[-\vareps^N,\vareps^N]} z\bone_{\{p>1,\La\notin\cals\}}\,\pi(\dd z)\right| \leq \calo\left((\vareps^N)^{1-q}\bone_{\{p>1,\La\notin\cals\}}\right)\;.
\end{align*}
For instance, if $\vareps^N=1/N^k$, then the sequence $(r^N_2)^{p^\ast} = \calo(N^{-k(p-q)})$ is summable for all $k>(p-q)^{-1}$. So in order to obtain a.s. convergence of $Y^N(t,x)\to Y(t,x)$, a sufficient condition is to choose the truncation rates $\vareps^N$ small enough. Similar conclusions are valid for the other two rates $r^N_1(T,K)$ and $r^N_3$.
\eexam

\section{Truncation via Series Representations}\label{s4}

From the viewpoint of simulation, the truncation of the small jumps as presented in the previous section, has two main drawbacks: first, it may not be so easy to simulate the jumps of the truncated L\'evy measure, i.e. from the distribution $\pi/\pi([-\vareps,\vareps]^\Comp)$, for a large number of times; second, the jumps have to be simulated all over again when one goes from step $N$ to step $N+1$. These two problems can be overcome by using series representations for the L\'evy basis. The idea, going back to \cite{Rosinski89,Rosinski90} and already applied to the simulation of L\'evy processes \cite{Rosinski01}, is to choose the jumps to be simulated in a random order. Instead of selecting the big jumps first and the smaller jumps later as in Sect.~\ref{s3}, we only choose the big jumps first more \emph{likely}. The details are as follows: we fix a finite time horizon $T\in\bbr_+$ and, recalling A7, a partition $(Q^i)_{i\in\bbn}$ of $\bbr^d$ into pairwise disjoint compact sets such that $U^N=\bigcup_{i=1}^N Q^i$. We now assume that the jump measure $\mu$ of $\La$ on the strip $[0,T]\times\bbr^d\times\bbr$ can be represented in the form
\begin{align} \mu(\dd t, \dd x, \dd z) &= \sum_{i=1}^\infty \mu_i(\dd t,\dd x,\dd z)\;,\nonumber \\
\mu_i(\dd t,\dd x,\dd z) &= \sum_{j=1}^\infty \delta_{(\tau^i_j,\xi^i_j,H(\Gamma^i_j,V^i_j))}(\dd t,\dd x,\dd z)\quad\text{a.s.}\;,\label{seriesrep} \end{align}
where $H\colon (0,\infty)\times \bbr \to \bbr$ is a measurable function, satisfying $H(\cdot,v)=-H(\cdot,-v)$ for all $v\in\bbr$ when $\La\in\cals$, 
and the random variables involved have the following properties for each $i\in\bbn$:
\bit
	\item $(\tau^i_{j}\colon j\in\bbn)$ and $(\xi^i_{j}\colon j\in\bbn)$ are i.i.d. sequences with uniform distribution on $[0,T]$ and $Q^i$, respectively.
	\item $(\Gamma^i_{j}\colon j\in\bbn)$ is a random walk whose increments are exponentially distributed with mean $1/T$.
	\item $(V^i_{j}\colon j\in\bbn)$ is an i.i.d. sequence with distribution $F$ on $\bbr$, which we should be able to simulate from. We assume that $F$ is symmetric when $\La\in\cals$. 
	\item The sequences $\tau^i$, $\xi^i$, $\Ga^i$ and $V^i$ are independent from each other.
	\item $(\tau^i,\xi^i,\Ga^i,V^i)$ is independent from $(\tau^k, \xi^k, \Ga^k, V^k\colon k\neq i)$.
\eit 
Because of \eqref{finvarsymm}, $\mu$ can always be written in the form \eqref{seriesrep} whenever the underlying stochastic basis is rich enough. We give three examples of such series representations.
\bexam\label{ex-sr} The proofs that the following choices are valid can be found in \cite[Sect.~3]{Rosinski01}, where also more examples are discussed.
\benu
	\item \emph{LePage's method}: we set $F:=(\delta_{-1}+\delta_1)/2$ and $H(r,\pm 1):=\pm\varrho^{\leftarrow}(r,\pm 1)$, where $\varrho^{\leftarrow}(r,\pm1)=\inf\{x\in(0,\infty)\colon \pi(\pm[x,\infty))<r\}$ for $r\in(0,\infty)$.
	\item \emph{Bondesson's method}: we assume that $\pi(A)=\int_0^\infty F(A/g(t))\,\dd t$ for $A\in\calb(\bbr^d)$ with some non-increasing $g\colon \bbr_+\to\bbr_+$. Then we define $H(r,v):=g(r)v$.
	\item \emph{Thinning method}: we choose $F$ in such a way that $Q$ is absolutely continuous with respect to $F$ with density $q$, and define $H(r,v):=v\bone_{\{q(v)\geq r\}}$.
\eenu
\eexam

Our approximation scheme is basically the same as in Sect.~\ref{s3}: we define $G^N$ by \eqref{GN} and $Y^N$ by \eqref{scheme1}, with the difference that $\La^N$ on $[0,T]\times\bbr^d$ is now defined as
\begin{align} \La^N(\dd t,\dd x)&:=\int_\bbr z\,\mu^N(\dd t,\dd x,\dd z)\;,\nonumber\\
 \mu^N(\dd t,\dd x,\dd z)&:=\sum_{i=1}^\infty \sum_{j\colon \Ga^i_{j}\leq N} \delta_{(\tau^i_{j},\xi^i_{j},H(\Gamma^i_{j},V^i_{j}))}(\dd t,\dd x,\dd z)\;. \label{LaN2} \end{align}
We can therefore rewrite $Y^N(t,x)$ for $(t,x)\in[0,T]\times\bbr^d$ as
\beq\label{scheme2var} Y^N(t,x) = Y_0(t,x) + \sum_{i=1}^N \sum_{j\colon \Ga^i_{j}\leq N} G(t,x;\tau^i_{j},\xi^i_{j})\si(Y^N(\tau^i_{j},\xi^i_{j})) H(\Ga^i_{j},V^i_{j})\bone_{\{\tau^i_{j}<t\}}\;. \eeq
This yields the following simulation algorithm:
\begin{alg}\label{alg2} Let $\calg$ be a finite grid in $[0,T]\times\bbr^d$ and $N\in\bbn$.
\benu \item For each $i=1,\ldots,N$ set $j:=1$ and repeat the following:
	\benu
		\item Draw $E^i_j$ from an exponential distribution with mean $1/T$.
		\item Define $\Ga^i_j:=\Ga^i_{j-1} + E^i_j$ ($\Ga^i_0:=0$).
		\item If $\Ga^i_j>N$, set $J_i:=j-1$ and leave the loop; otherwise set $j:=j+1$.
	\eenu
\item For each $i=1,\ldots,N$ and $j=1,\ldots,J_i$ simulate independently
	\benu
		\item a pair $(\tau^i_j,\xi^i_j)$ with uniform distribution on $[0,T]\times Q^i$;
		\item a random variable $V^i_j$ with distribution $F$;
		\item the random variable $Y_0(\tau^i_j,\xi^i_j)$.
	\eenu
\item Sort the sequence $(\tau^i_j\colon i=1,\ldots,N, j=1,\ldots,J_i)$ in increasing order, yielding sequences $(\tau_i,\xi_i,\Ga_i,V_i\colon i=1,\ldots, \sum_{j=1}^N J_j)$. Now define
\[ Y^N(\tau_i,\xi_i):=Y_0(\tau_i,\xi_i) + \sum_{j=1}^{i-1} G(\tau_i,\xi_i;\tau_j,\xi_j)\si(Y^N(\tau_j,\xi_j))H(\Ga_j,V_j)\;. \]
\item For each $(t,x)\in\calg$ simulate $Y_0(t,x)$ and define $Y(t,x)$ by \eqref{scheme2var}.
\eenu
\end{alg}
We can now prove a convergence theorem for $Y^N$ to $Y$, similar to Thm.~\ref{conv1}. Define
\begin{align} r^N_2&:= \left(\int_N^\infty \int_\bbr |H(r,v)|^p\,F(\dd v)\,\dd r\right)^{1/{p^\ast}}\;,\nonumber\\
 r^N_3 &:= \left|\int_N^\infty \int_\bbr H(r,v)\bone_{\{p>1,\La\notin\cals\}} \,F(\dd v)\,\dd r\right|\;. \label{rN-2}\end{align}
\bthe\label{conv2} Grant assumptions A1--A7 under which the SVE \eqref{stochint} has a unique solution in $B^p_\loc$. Further suppose that the jump measure $\mu$ of $\La$ has a representation in form of \eqref{seriesrep}. Then $Y^N$ as defined in \eqref{scheme2var} belongs to $B^p_\loc$ for all $N\in\bbn$, and for all $T\in\bbr_+$ and compact sets $K\subseteq\bbr^d$ there exists a constant $C(T)\in\bbr_+$ independent of $N$ and $K$ such that
\beq\label{convtrunc2} \sup_{(t,x)\in[0,T]\times K} \|Y(t,x)-Y^N(t,x)\|_{L^p} \leq C(T) (r^N_1(T,K) + r^N_2 + r^N_3)\;. \eeq
If $\sum_{n=1}^\infty (r^N_1(T,K)+r^N_2 +r^N_3)^{p^\ast} < \infty$, then we also have for all $(t,x)\in[0,T]\times K$ that $Y^N(t,x)\to Y(t,x)$ a.s. as $N\to\infty$.
\ethe
\bproof We start with some preliminaries. It follows from \eqref{seriesrep} and \cite[Prop.~2.1]{Rosinski01} that on $[0,T]\times\bbr^d\times\bbr$ we have $\nu=\bar\nu \circ h^{-1}$ where $\bar\nu(\dd t,\dd x,\dd r,\dd v) = \dd t\,\dd x\,\dd r\,F(\dd v)$ and $h(t,x,r,v)=(t,x,H(r,v))$. Therefore, conditions \eqref{finvarsymm} and \eqref{intp} imply that
\[ \int_0^\infty \int_\bbr |H(r,v)|^p\,F(\dd v)\,\dd r = \int_\bbr |z|^p\,\pi(\dd z) < \infty\;, \]
\[ \text{and}\quad \int_0^\infty \int_\bbr |H(r,v)|\bone_{\{p>1,\La\notin\cals\}}\,F(\dd v)\,\dd r = \int_\bbr |z|\bone_{\{p>1,\La\notin\cals\}}\,\pi(\dd z)<\infty\;. \]
Consequently, $r^N_2$ and $r^N_3$ are well defined and converge to $0$ when $N\to\infty$. Similarly, the compensator $\nu_N$ of the measure $\mu-\mu^N$ is given by $\nu_N(\dd t,\dd x,\dd z)=\dd t\,\dd x\,\pi_N(\dd z)$, where $\pi_N=(\Leb\otimes F)\circ H_N^{-1}$ and $H_N(r,v)=H(r,v)\bone_{(N,\infty)}(r)$.

For the actual proof of Thm.~\ref{conv2} one can basically follow the proof of Thm.~\ref{conv1}. Only the estimation of $I^N_3(t,x)$ as defined in \eqref{INs} is different, which we shall carry out now. In the case of $p>1$, we again use the Burkholder-Davis-Gundy inequality and obtain for $(t,x)\in[0,T]\times\bbr^d$
\begin{align*}
\|I^N_3(t,x)\|_{L^p} &=  \left\|\int_0^t\int_{\bbr^d}\int_\bbr G^N(t,x;s,y) \si(Y^N(s,y))z\,(\mu_N-\nu_N)(\dd s,\dd y,\dd z)\right\|_{L^p}\nonumber \\
 &\quad+\left\|\int_0^t\int_{\bbr^d} \int_\bbr G^N(t,x;s,y) \si(Y^N(s,y)) z\bone_{\{\La\notin\cals\}} \,\nu_N(\dd s,\dd y,\dd z)\right\|_{L^p}\\
 &\leq C(T)\Bigg(\left(\int_0^t\int_{\bbr^d} \int_N^\infty\int_\bbr |G(t,x;s,y)H(r,v)|^p \,F(\dd v)\,\dd r\,\dd(s,y)\right)^{1/p} \\
 &\quad+ \int_0^t\int_{\bbr^d} |G(t,x;s,y)| \left|\int_N^\infty \int_\bbr H(r,v)\bone_{\{\La\notin\cals\}}\,F(\dd v)\,\dd r\right|\,\dd(s,y)\Bigg)\\
 &\leq C(T) (r^N_2+r^N_3)\;.
\end{align*}
The case $p\in(0,1]$ is treated analogously. One only needs to replace $\mu_N-\nu_N$ by $\mu_N$ and estimate via Jensen's inequality. \qed
\eproof

\bexam[Continuation of Ex.~\ref{ex-sr}] We calculate the rates $r^N_2$ and $r^N_3$ from \eqref{rN-2} for the series representations given in Ex.~\ref{ex-sr}. We assume that $p,q\in(0,2]$ with $q<p$ are chosen such that $\int_{[-1,1]} |z|^q\,\pi(\dd z)<\infty$, and $q<1$ if $\La\in\calv_0$. For all three examples we use the fact that $\pi=(\Leb\otimes F)\circ H^{-1}$ and that $r>N$ implies $|H(r,v)|\leq |H(N,v)|$ for all $v\in\bbr$.
\benu
	\item \emph{LePage's method}: We have
	\begin{align*} (r^N_2)^{p^\ast}&=\int_N^\infty \frac{|H(r,1)|^p + |H(r,-1)|^p}{2}\,\dd r\leq \frac{1}{2} \int_{[H(N,-1),H(N,1)]} |z|^p\,\pi(\dd z)\\
	&\leq \frac{1}{2} \int_{[H(1,-1),H(1,1)]} |z|^q\,\pi(\dd z) (|H(N,-1)|\vee|H(N,1)|)^{p-q}\;,\end{align*}
	and therefore
	\begin{align*}
	r^N_2&= \calo\left((\varrho^\leftarrow(N,1)\vee \varrho^\leftarrow(N,-1))^{(p-q)/p^\ast}\right)\;,\\
	r^N_3 &= \calo\left((\varrho^\leftarrow(N,1)\vee \varrho^\leftarrow(N,-1))^{1-q}\bone_{\{p>1,\La\notin\cals\}}\right)\;.\end{align*}
	\item \emph{Bondesson's method}: Since $H(r,v)=g(r)v$ and $g$ is non-increasing, we obtain
	\[(r^N_2)^{p^\ast}=\int_N^\infty \int_\bbr |g(r)v|^p\,F(\dd v)\,\dd r \leq (g(N))^{p-q}\int_0^\infty g^q(r)\,\dd r\int_\bbr |v|^p\,F(\dd v)\;,\]
	and consequently
	\[ r^N_2=\calo\left(g(N)^{(p-q)/{p^\ast}}\right)\;,\quad r^N_3 = \calo\left(g(N)^{1-q}\bone_{\{p>1,\La\notin\cals\}}\right)\;.\]
	\item \emph{Thinning method}: Here we have
	\begin{align*} r^N_2 &=\left(\int_\bbr\int_N^{q(v)\vee N}  |v|^p \,\dd r\,F(\dd v)\right)^{1/p^\ast} = \left( \int_\bbr |v|^p \frac{(q(v)-N)\vee 0}{q(v)}\,\pi(\dd v) \right)^{1/p^\ast}\\
	&\leq \left(\int_\bbr |z|^p \bone_{\{q(v)\geq N\}}\,\pi(\dd z)\right)^{1/p^\ast}\;,\\
	r^N_3 &\leq \int_\bbr |z| \bone_{\{q(v)\geq N\}} \bone_{\{p>1,\La\notin\cals\}}\,\pi(\dd z)\;.
	\end{align*}
	In most situations, there exist $(\vareps^N)_{N\in\bbn}\subseteq \bbr_+$ with $\vareps^N\to0$ as $N\to\infty$ such that $\{q(v)\geq N\} \subseteq [-\vareps^N,\vareps^N]$. In this case, one can apply the estimates in Ex.~\ref{ex1}.
\eenu\eexam

\section{Simulation Study}\label{s5}
\begin{figure}[ht!]
  \centering
  \includegraphics[width=\linewidth]{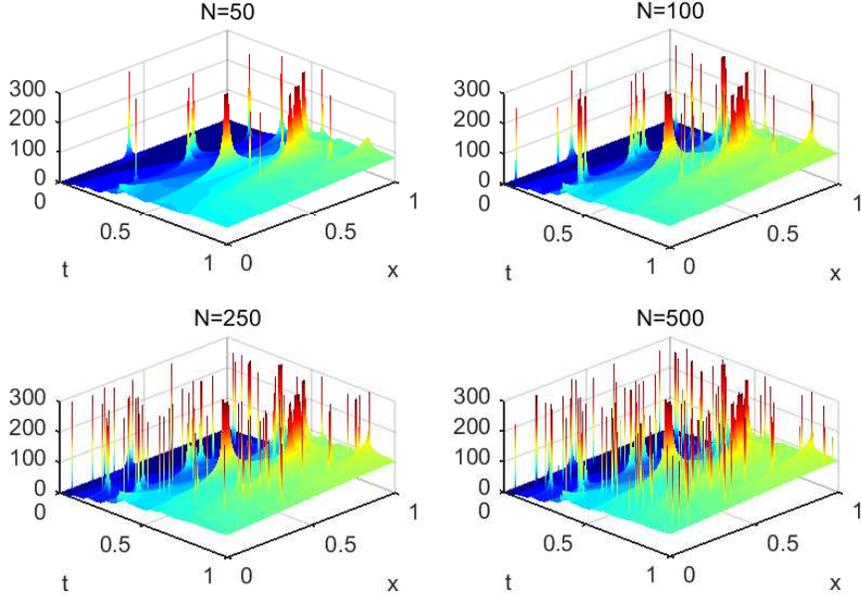}
  \caption{Successive approximations of $Y$ as given in \eqref{simple} via Bondesson's method in dimension $1$ for $(t,x)\in[0,1]\times[0,1]$ with $N\in\{50,100,250,500\}$ jumps in the region $[0,1]\times[-1,2]$}\label{bond-dim1}
\end{figure}
In this section we visualize the sample path behaviour of the stochastic heat equation from Ex.~\ref{ex1} via a simulation study, using MATLAB programs from \cite{Chen14}. We take $\La$ to be a L\'evy basis without drift, whose L\'evy measure $\pi$ is that of a gamma process, i.e. 
\[\pi(\dd z)=\gamma z^{-1}\exp{(-\lambda z)}\bone_{\{z>0\}}\,\dd z\]
with two parameters $\ga,\la>0$. In the figures below their values are always $\ga=10$ and $\la=0.1$. Furthermore, we set $Y_0\equiv0$ and $\si\equiv1$.
Especially the latter choice simplifies the subsequent discussion a lot, but none of the issues we address below relies on this assumption. Thus, the process we would like to simulate is
\beq\label{simple} Y(t,x)=\int_0^t\int_{\bbr^d} g(t-s,x-y)\,\La(\dd s,\dd y)\;,\quad(t,x)\in\bbr_+\times\bbr^d\;, \eeq
with $g$ being the heat kernel given in \eqref{hk}.
In order to understand the path properties of $Y$, it is important to notice that $g$ is smooth on the whole $\bbr_+\times\bbr^d$ except at the origin where it explodes. More precisely, for every $t\in(0,\infty)$ the function $x\mapsto g(t,x)$ is the Gaussian density with mean $0$ and variance $2t$, which is smooth and assumes its maximum at $0$. Also, for every $x\neq0$, the function $t\mapsto g(t,x)$ is smooth (also at $t=0$), with maximum at $t=\|x\|^2/(2d)$. However, if $x=0$, then $g(t,0)=(4\pi t)^{-d/2}$ has a singularity at $t=0$. 

These analytical properties have direct consequences on the sample paths of $Y$. When $\La$ is of compound Poisson type, that is, has only finitely many atoms on compact sets, it can be readily seen from \eqref{simple} that the evolution of $Y$ after a jump $J$ at $(\tau,\xi)$ follows the shape of the heat kernel until a next jump arrives. In particular, for $x=\xi$, $Y(t,x)$ jumps to infinity at $\tau$, and decays in $t$ like $J(4\pi (t-\tau))^{-d/2}$ afterwards. But for every $x\neq\xi$, the evolution $t\mapsto Y(t,x)$ is \emph{smooth} at $t=\tau$. In fact, it first starts to increase until $t=\tau+\|x-\xi\|^2/(2d)$ and then decays again. As a consequence, in space dimension $1$, the space--time plot of $Y$ shows a basically smoothly evolving path, only interrupted with slim poles at the jump locations of $\La$; see the case $N=50$ in Fig.~\ref{bond-dim1}.
However, when $\La$ has infinite activity, that is, has infinitely many jumps on any non-empty open set, then it is known from \cite[Thm.~4]{Rosinski89} that on any such set $Y$ is unbounded, at least with positive probability. Therefore, the space--time plots of the approximations of $Y$ with finitely many jumps must be treated with caution: in the limiting situation, no smooth area exists any more, but there will be a dense subset of singularities on the plane, which is in line with Fig.~\ref{bond-dim1}. 

\begin{figure}[t!]
  \centering
  \includegraphics[width=\linewidth]{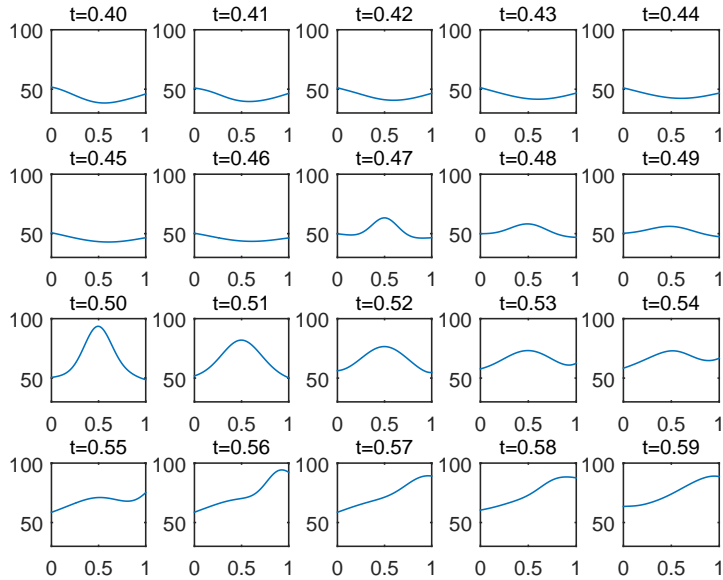}
  \caption{Several $t$-sections of the realization of $Y$ shown in Fig.~\ref{bond-dim1} with $N=500$}\label{section-dim1}
\end{figure}

\begin{figure}
  \centering
  \includegraphics[width=\linewidth]{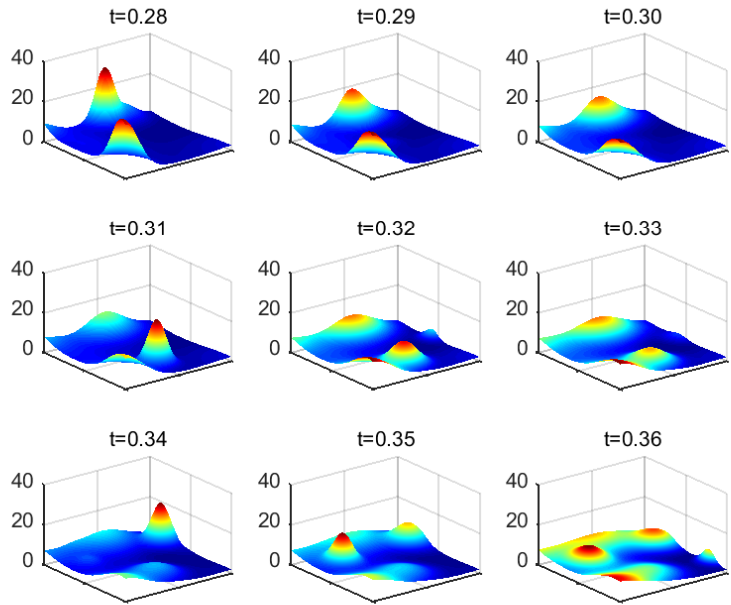}
    \vspace{-3\baselineskip}
  \caption{Several $t$-sections in the region $[-1,1]^2$ of a realization of $Y$ in dimension $2$ by Bondesson's method with $N=500$ jumps within $[0,1]\times[-2,2]^2$}\label{section-dim2}
\end{figure}

\begin{figure}
\vspace{2\baselineskip}
  \centering
  \includegraphics[width=\linewidth]{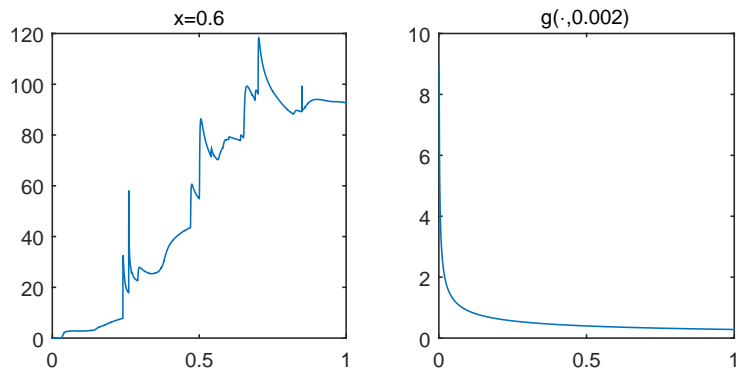}
  \caption{The $x$-section of the realization of $Y$ as in Fig.~\ref{bond-dim1} with $N=500$ at $x=0.6$ and the heat kernel $g(\cdot,x)$ at $x=0.002$}\label{x-section-dim1}
  \end{figure}

Another interesting observation, however, is the following: if we consider a countable number of $x$- or $t$-sections of $Y$ (for $x\in\bbr^d$, the $x$-section of $Y$ is given by the function $t\mapsto Y(t,x)$; for $t\in\bbr_+$, the $t$-section of $Y$ is the function $x\mapsto Y(t,x)$), then it is shown in \cite[Sect.~2]{SLB98} that these are continuous with probability one. Intuitively, this is possible because a.s. the sections never hit a jump (although they are arbitrarily close). For instance, Figs.~\ref{section-dim1} and \ref{section-dim2} show $t$-sections of a realization of \eqref{simple} in one, respectively two space dimensions. So as long as we only take countably many ``measurements'', we do not observe the space--time singularities of $Y$ but only its relatively regular sections. In theory, this also includes the $x$-sections of the process $Y$. But if we plot them for one space dimension as in Fig.~\ref{x-section-dim1}, one would conjecture from the simulation that they exhibit jumps in time. However, this is \emph{not} true: the jump-like appearance of the $x$-sections are due to the fact that $g(\cdot,x)$ resembles a discontinuous function at $t=0$ for small $x$. Of course, it follows right from the definition \eqref{hk} that all $x$-sections of $g$ are smooth everywhere. 

\section{Conclusion and Outlook}\label{s6}
In Sects.~\ref{s3} and \ref{s4} we have presented two simulation algorithms for the SVE \eqref{stochint}: Algorithms~\ref{alg1} and \ref{alg2}. In Thms.~\ref{conv1} and \ref{conv2} we have determined the rate of convergence of the approximations $Y^N$ to $Y$ in the $L^p$-sense. If these rates are small enough, we have also proved a.s. convergence. Although the theoretical analysis of both schemes lead to quite similar results regarding their convergence behaviour, there are important differences which will decide on whether the one or the other method is preferable in concrete situations. For the first method of truncating the small jumps to work, one must be able to efficiently simulate from the truncated L\'evy measure $\pi/\pi([-\vareps,\vareps]^\Comp)$ for small $\vareps$. 
For the second method, which relies on series representations, the main challenge is to choose $H$ and $F$ in a way such that $H$ is explicitly known \emph{and} $F$ can be easily simulated from. For instance, if one uses LePage's method (see Ex.~\ref{ex-sr}), then $F=(\delta_{-1}+\delta_1)/2$ is easily simulated, but for $H$, which is given by the generalized inverse tails of the underlying L\'evy measure, maybe no tractable expression exists.

Finally, let us comment on further generalizations of the our results. Throughout this paper, we have assumed that the driving noise $\La$ is a homogeneous L\'evy basis, i.e. satisfies \eqref{invar}. In fact, we have introduced this condition only for the sake of simplicity: with a straightforward adjustment, all results obtained in this paper also hold for time- and space-varying (but deterministic) characteristics. Another issue is the finite time perspective which we have taken up for our analysis. An interesting question would be under which conditions \eqref{stochint} has a stationary solution, and in this case, whether one can simulate from it. Sufficient conditions for the existence and uniqueness of stationary solutions to \eqref{stochint} are determined in \cite[Thm.~4.8]{Chong14b}. Under these conditions, the methods used to derive Thms.~\ref{conv1} and \ref{conv2} can indeed be extended to the case of infinite time horizon. We leave the details to the reader at this point. 

At last, also the hypothesis that $\La$ is of pure-jump type can be weakened. If $\La$ has an additional drift (including the case where $\La$ has locally infinite variation and is not symmetric) but still no Gaussian part, the approximations $Y^N$ in \eqref{scheme1} or \eqref{scheme2var} will contain a further term that is a Volterra integral with respect to the Lebesgue measure. So each time in between two simulated jumps, a deterministic Volterra equation has to be solved numerically, which boils down to a deterministic PDE in the case where $G$ comes from an SPDE. For this subject, there exists a huge literature, which is, of course, also very different to the stochastic case as considered above. If $\La$ also contains a Gaussian part, then one has to apply techniques from the papers cited in Sect.~\ref{s1} and ours simultaneously. We content ourselves with referring to \cite{Zhang08}, who numerically analyzes a Volterra equation driven by a drift plus a Brownian motion. Finally, let us remark that if $p=2$ (in particular, $G$ must be square-integrable), it is possible for some L\'evy bases to improve the results of Sect.~\ref{s3} if we do not neglect the small jumps completely but approximate them via a Gaussian noise with the same variance, cf. \cite{Asmussen01} in the case of L\'evy processes.

\begin{acknowledgement}
We take pleasure in thanking Jean Jacod for his valuable advice on this subject. The second author acknowledges support from the Studienstiftung des deutschen Volkes and the graduate programme TopMath at Technische Universit\"at M\"unchen.
\end{acknowledgement}

\bibliographystyle{spmpsci}
\bibliography{bib-CCK}
\end{document}